\documentclass[12pt,a4paper,reqno,twoside]{amsart}

\usepackage[english]{babel}
\usepackage{stmaryrd}
\usepackage{dsfont}
\usepackage{bbm}
\usepackage[symbol*,ragged]{footmisc}
\usepackage[colorlinks,linkcolor=red,anchorcolor=blue,citecolor=blue,urlcolor=blue]{hyperref}
\usepackage{color,xcolor}

\usepackage{geometry}
\usepackage{amssymb}
\usepackage{amsmath}
\usepackage{mathrsfs}
\usepackage{amsfonts}
\usepackage{epsfig}

\usepackage{amsthm}
\usepackage{amsxtra}
\usepackage{bbding}
\usepackage{epsfig}
\usepackage{graphicx}
\usepackage{latexsym}
\usepackage{mathbbol}
\usepackage{bbold}

\usepackage{pifont}
\usepackage{wasysym}
\usepackage{skull}
\usepackage{float}

\DeclareSymbolFontAlphabet{\mathbb}{AMSb}
\DeclareSymbolFontAlphabet{\mathbbol}{bbold}

\usepackage{amscd}
\usepackage[all]{xy}
\allowdisplaybreaks[4]
\usepackage{setspace}

\theoremstyle{plain}

\newtheorem*{corollary*}{\normalfont\scshape Corollary}

\theoremstyle{remark}
\newtheorem*{remark*}{\normalfont\scshape Remark}

\numberwithin{equation}{section}
\addtocounter{footnote}{1}

\renewcommand{\footnoterule}{
  \kern -3pt
  \hrule width 2.5in height 0.4pt
  \kern 3pt
}

\makeatletter
\@ifundefined{MakeUppercase}{}{}
\makeatother

\begin{document}
	
\title[A survey on moments of a class of error terms ]
{A survey on moments of a class of error terms  } 

\author{Wenguang Zhai \\ \\ published in {\it Adv. Math.} (China) {\bf 55} (2026), 481-504}

\address{[Wenguang Zhai]  Department of Mathematics, China University of Mining and Technology, Beijing 100083, People's Republic of China}

\email{zhaiwg@hotmail.com}

\date{}

	

\thanks{This work is  partially supported  by
the National Natural Science Foundation of China(Grant No. 12471009).}

\begin{abstract}
  This is a survey article on the moments of a class of error terms in analytic number theory. We begin by reviewing  the moments of the Dirichlet divisor problem,  and then review the hybrid moments of error terms in the divisor problems. We give two new results on the joint distribution of sign changes of $\Delta_2(x)$ and $\Delta_3(x). $ In the last section, we present some conjectures on the moment  results in the well-known Selberg  class.
\end{abstract}

\maketitle

\section{Introduction}

\subsection{The mean value problem  in number theory}\

Let $f: \mathbb{N}\rightarrow   \mathbb{C} $ be any arithmetic function. Without loss of generality, we suppose
that $f(n)\ll n^\varepsilon$ holds for any $\varepsilon>0$. It is a classical problem in the analytic number theory to study the asymptotic behaviour of the mean value
\begin{equation*}
S_f(x):=\sum_{n\leq x}f(n),\ \ x\rightarrow \infty.
\end{equation*}

Let
\begin{equation*}
\mathfrak{L}_f(s):=\sum_{n=1}^\infty\frac{f(n)}{n^s}\ \ (\Re(s)>1).
\end{equation*}

Suppose $\mathfrak{L}_f(s)$  can be continued meromorphically to the whole complex plane $\mathbb{C}$ and has a pole of degree $d_f>0$ at $s = 1$. The expected asymptotic formula of $S_f(x)$ is of the form
\begin{equation*}
S_f(x)=M_f(x)+E_f(x),
\end{equation*}
where
\begin{equation*}
M_f(x)= Res_{s=1}\mathfrak{L}_f(s) x^ss^{-1}
\end{equation*}
is the main term and $E_f(x)$ is the error term.

The history of this kind of problem goes back to Gauss and Dirichlet. Gauss first
proved the asymptotic formula
\begin{equation}\label{yi.1}
 \sum_{n\leq x}r(n)=\pi x+P(x), \ \ P(x)=O(x^{1/2}),
\end{equation}
where $r(n)$ denotes the number of ways $n$ can be written as a sum of two squares.
Dirichlet first proved that
\begin{equation}\label{yi.5}
 \sum_{n\leq x}d(n)=  x\log x+(2\gamma-1)x+\Delta(x), \ \ \Delta(x)=O(x^{1/2}),
\end{equation}
where $d(n)$ denotes the number of ways $n$ can be written as a product of two natural
numbers and $\gamma$ is the Euler constant.

An important problem in analytic number theory is to study various properties
of $E_f(x)$ for different $f^{'}$s. For example, one may investigate its upper bound, its
Omega estimate, its moments, its sign changes and its probability distribution, etc.
For each of the above aspects, there are abundant results in analytic number theory.

In this paper, we will give a survey of moments results and applications of
$E_f(x)$
for an important class of arithmetic functions in number theory, whose Dirichlet
series form the well-known Selberg class.

\subsection{The Selberg class}\

The Selberg class $\mathcal{S}$(see for
example \cite{CG, Kac, KP, KP2, Ka, Se}) consists of all non-vanishing Dirichlet series
\begin{equation}\label{yi.3}
\mathfrak{L}_f(s) :=\sum_{n=1}^\infty\frac{f(n)}{n^s},
\end{equation}
which satisfies the following hypotheses:

{\bf I. Ramanujan's conjecture}: $f(n)\ll n^\varepsilon$ for any $\varepsilon>0$;

{\bf  II. Analytic continuation} : there exists a non-negative integer
$m_{\mathfrak{L}_f}$ such that $(s-1)^{m_{\mathfrak{L}_f}}\mathfrak{L}_f(s)$ is an entire function of finite order;

{\bf III. Functional equation}: $\mathfrak{L}_f(s)$ satisfies a functional equation of type  $$\Lambda_{\mathfrak{L}_f}(s)=\omega \overline{\Lambda_{\mathfrak{L}_f}(1-\overline{s})},$$
where
$$ \Lambda_{\mathfrak{L}_f}(s)=\mathfrak{L}_f(s)Q^s\prod_{j=1}^{\mathcal{R}}\Gamma(\lambda_j s+\mu_j)$$
with positive real numbers $Q, \lambda_j$ and complex numbers $\mu_j, \omega$ with $\Re{\mu_j}\geq 0$ and
$|\omega|=1$.

{\bf IV. Euler product}: $\mathfrak{L}_f(s)$ satisfies
$$\mathfrak{L}_f(s)=\prod_p \exp\left(\sum_{n\geq 1}\frac{b(p^n)}{p^{ns}}\right)$$
with suitable coefficients $b(p^n)$ satisfying $b(p^n)\ll p^{nc}$ for some
$c < 1/2.$

The number $d = 2\sum_{j}\lambda_j$ is called the degree of $\mathfrak{L}_f(s)$.
It is conjectured that all functions in $\mathcal{S}$ have integral degrees.
It is well-known that a product of a function $\mathfrak{L}_{f_1}(s)$ in $\mathcal{S}$ of degree $k_1$ and a function $\mathfrak{L}_{f_2}(s)$ in $\mathcal{S}$ of degree $k_2$ is a function in $\mathcal{S}$ of degree $k_1k_2.$

Many well-known functions are contained in the Selberg class $\mathcal{S}$. We recall some examples.

(1) The famous Riemann zeta-function $\zeta(s)$ is a function in $\mathcal{S}$ of degree $1$. Actually for any fixed integer $k\geq 1$, $ \zeta^k(s)$ is a function in $\mathcal{S}$ of degree $k$.

(2) Let $g(z)$ be a holomorphic cusp form with respect to $SL_2(\mathbb{Z})$, then the auto-morphic $L$-function attached to $g(z)$ is a function in $\mathcal{S}$ of degree $2.$

(3) Let $F$ be an algebraic number field of degree $k.$ Then the Dedekind
zeta-function $\zeta_F (s)$ is a function in $\mathcal{S}$ of degree $k.$

{\bf Notation.}  Throughout this paper,  ${\Bbb N}$ denotes the set of positive integers, ${\Bbb R} $ denotes the set of real numbers  and  ${\Bbb C}$ denotes the set of  all complex numbers.
For any real number $t,$ let  $\Vert t\Vert $ denote the distance between $t$ and the nearest integer. Let $\zeta(s)$ denote the Riemann zeta function and $\varepsilon$ denote a small positive constant. For $n\in {\Bbb N},$ let $d_k(n)$ denote the $k$-fold divisor function and $d_2(n)=d(n).$
  The symbols  $f=O(g)$ and $f\ll g$ mean
there exists a positive constant  $c>0 $ such that
$ |f|\leq cg. $

\section{ Moments of $\Delta(x)$}

One of the main tasks of the Dirichlet divisor problem is to use the method of
exponential sums to reduce the exponent $1/2$ in the error term of (\ref{yi.5}). The study of the
Dirichlet divisor problem promoted the development of the methods of exponential
sums in analytic number theory. The Dirichlet's exponent $1/2$ was improved by
many authors. The  latest result reads(see, Huxley \cite{Hux})
\begin{equation}\label{er.1}
\Delta(x)\ll   x^{131/416}(\log x)^{26957/8320}.
\end{equation}
It is conjectured that
\begin{equation}
\Delta(x)\ll x^{1/4+\varepsilon},
\end{equation}
which is supported by the moment results of $\Delta(x)$.

\subsection{Moments of $\Delta(x)$}\

Vorono\"{\i} in \cite{Vo} proved
\begin{equation}\label{er.3}
\int_1^T\Delta(x)dx = \frac{T}{4}+\frac{T^{3/4}}{2\sqrt 2 \pi^2}
\sum_{n=1}^\infty\frac{d(n)}{n^{5/4}}
\sin(4\pi\sqrt {nT}-\frac{\pi}{4})+O(T^{1/4 }).
\end{equation}

In 1922 Cram\'{e}r \cite{Cra} proved the mean square result
\begin{equation}\label{er.4}
\int_1^T\Delta^2(x)dx =C_2^{(2)}T^{3/2}+O(T^{5/4+\varepsilon}),\ \ C_2^{(2)}=\frac{\zeta^4(3/2)}{6\pi^2 \zeta(3)}.
\end{equation}

The proof of (\ref{er.4}) relies on the following well-known Vorono\"{\i}'s formula(see, for
example, \cite{Iv} ), which is the basis for the study of the $\Delta(x).$

{\bf Lemma 2.1. } {\it Suppose $1\ll N\ll x.$ Then
\begin{equation*}
\Delta(x)=\frac{x^{1/4}}{\sqrt 2 \pi} \sum_{n\leq N}\frac{d(n)}{n^{3/4}}
\cos(4\pi\sqrt {nx}-\frac{\pi}{4}) +O(x^{1/2+\varepsilon}N^{-1/2} ).
\end{equation*} }

Let $F_2(T)$ denote the error term in (\ref{er.4}). In 1956,  Tong \cite{To2} proved that
$$F_2(T) = O(T \log^5 T).$$   Preissman \cite{Pr} proved that
\begin{equation}\label{er.5}
 F_2(T)= O(T \log^4 T).
\end{equation}
In \cite{LT2} Lau and Tsang proved
\begin{equation}\label{er.6}
F_2(T) = O(T \log^3 T\log\log T),
\end{equation}
which is the best upper bound of $F_2(T)$ so far.

The proof of (\ref{er.5}) and the proof of (\ref{er.6}) depend on the Vorono\"{\i}'s formula of the  following form (see, Meurman \cite{Meu}). Lemma 2.2 also plays an important role in  higher moments of $\Delta(x).$

{\bf Lemma 2.2. } {\it Suppose $x\ll N\ll x^\Theta$ for some $\Theta>1.$ Then
\begin{equation*}
\Delta(x)=\frac{x^{1/4}}{\sqrt 2 \pi} \sum_{n\leq N}\frac{d(n)}{n^{3/4}}
\cos(4\pi\sqrt {nx}-\frac{\pi}{4}) +R(x,N),
\end{equation*}
where
\begin{eqnarray*}
R(x,N)=\left\{\begin{array}{ll}
x^\varepsilon,&\mbox{if $\Vert x\Vert\ll x^{5/2}N^{-1/2},$}\\
x^{-1/4},&\mbox{if $\Vert x\Vert\gg x^{5/2}N^{-1/2}.$}
\end{array}\right.
\end{eqnarray*} }

Let $A >0$ be a fixed real number.   Ivi\'{c} \cite{Iv2} proved that the
estimate
\begin{eqnarray}\label{er.7}
 \int_1^T|\Delta(x)|^Adx\ll T^{1+A/4+\varepsilon}
\end{eqnarray}
holds for $0<A\leq  8.75$(See also Theorem 13.9 of Ivi\'{c} \cite{Iv}). Ivi\'{c} first proved a large value estimate of $\Delta(x)$ (see Theorem 1 of Ivi\'{c} \cite{Iv2} or Theorem 13.8 of Ivi\'{c} \cite{Iv}), then from which he deduced the estimate (\ref{er.7}). If we insert Huxley's estimate (\ref{er.1}) into Ivi\'{c}'s machinery, we see that (\ref{er.7}) holds for $0 < A \leq 262/27 = 9.\overline{703}. $

Using Lemma 2.1, in 1990 Tsang \cite{Ts} first proved that the asymptotic formulas
\begin{eqnarray}\label{er.8}
 \int_1^T \Delta^3(x) dx=C_2^{(3)} T^{7/4}+O(T^{7/4-1/14+\varepsilon})
\end{eqnarray}
and
\begin{eqnarray}
 \int_1^T \Delta^4(x) dx=C_2^{(4)} T^{2}+O(T^{2-1/23+\varepsilon})
\end{eqnarray}
hold, where $C_2^{(3)}$ and $C_2^{(4)}$ are two positive real numbers. Here $C_2^{(3)} > 0$ shows that the values of $\Delta(x)$ have a substantial bias towards the positive side.

Heath-Brown \cite{HB} proved that for any fixed $0 < A \leq 9$ the mean value
\begin{eqnarray}\label{er.10}
 T^{-1-A}\int_1^T |\Delta(x)|^A dx
\end{eqnarray}
converges to a finite limit as $T$ tends to infinity. Moreover, the same is true for the odd moments
\begin{eqnarray}
 T^{-1-k}\int_1^T \Delta^k(x) dx
\end{eqnarray}
for $k = 3, 5, 7, 9.$

Following Tsang's approach, Zhai \cite{Zw1} proved that
\begin{eqnarray}
&&\int_1^T \Delta^3(x) dx=C_2^{(3)} T^{7/4}+O(T^{7/4-1/4+\varepsilon}),\\
&&\int_1^T \Delta^4(x) dx=C_2^{(4)}T^{2}+O(T^{2-2/41+\varepsilon}),\nonumber\\
&&\int_1^T \Delta^5(x) dx=C_2^{(5)}T^{9/4}+O(T^{9/4- 5/816+\varepsilon}),\nonumber
\end{eqnarray}
where $C_2^{(5)}$ is a constant.

Ivi\'c and Sargos \cite{IS} proved that
\begin{eqnarray}\label{er.13}
&&\int_1^T \Delta^3(x) dx=C_2^{(3)} T^{7/4}+O(T^{7/4-7/20+\varepsilon}),\nonumber\\
&&\int_1^T \Delta^4(x) dx=C_2^{(4)} T^{2}+O(T^{2-1/12+\varepsilon}).
\end{eqnarray}

Let $k \geq 3$ be a fixed integer. Zhai \cite{Zw2} used a unified approach to study the higher moments of $\Delta(x)$. He proved that if $ A_0 > 9$ is a fixed real number such that the estimate   (\ref{er.7}) holds, then for any integer $3 \leq k < A_0$,  we can get an asymptotic formula  of $\int_1^T \Delta^k(x)dx.$ Suppose $T \leq x \leq 2T$ and
$y = T^\delta$ where $\delta > 0$ is a small parameter  to be determined later. Write (by Lemma 2.1)
\begin{equation*}
\Delta_1(x, y)=\frac{x^{1/4}}{\sqrt 2 \pi} \sum_{n\leq y}\frac{d(n)}{n^{3/4}}
\cos(4\pi\sqrt {nx}-\frac{\pi}{4}) ,\ \
\Delta_2(x, y)=\Delta(x)-\Delta_1(x, y).
\end{equation*}
Using Ivi\'c's large-value approach to $\Delta_2(x; y)$ directly we see
 $\Delta_2(x; y)$ is small on average. Since $\delta> 0$ is small, it is easy to get an asymptotic formula of $\int_T^{2T} \Delta^k_1(x; y)dx.$ Hence we can get an asymptotic formula of $\int_T^{2T} \Delta^k(x)dx$  by choosing a suitable $\delta$ to
 balance the contributions of $\Delta_1(x, y)$ and  $\Delta_2(x, y).$
Using this idea, Zhai \cite{Zw2} proved the asymptotic formula
\begin{eqnarray}\label{er.14}
 \int_1^T \Delta^k(x) dx=C_2^{(k)} T^{1+k/4}+O(T^{1+k/4- \delta_k+\varepsilon})
 \ \ (k=3, 4, 5, 6, 7, 8, 9),
\end{eqnarray}
where $C_2^{(k)}$ is a constant,  $\delta_5 = 1/64, \ \delta_6 = 35/4742,\
 \delta_7 = 17/6312,\ \delta_8 = 8/9433,\  \delta_9 = 13/75216.$ Furthermore, Zhai
 also pointed out that $C_2^{(k)}>0$ for $k=5, 7, 9.$

The values of $\delta_k\ (4\leq k\leq 9)$ were further studied by many authors. For example,  Zhai \cite{Zw3} proved   $\delta_4=3/28$ and   Kong \cite{kon} proved  $\delta_4=1/8.$   For better values of $\delta_k\ (5\leq k\leq 9),$ see
Zhang and Zhai \cite{ZZ} for $k=5,$ see Wang \cite{Wj} for $k=6,$
 see \cite{CTZ} for $k=7, 8, 9.$ See also Wang \cite{Wj} for $k=6 $
 and Li \cite{Li} for $k=7.$

For any integer $k\geq 2, $ let
$$F_k(T):=\int_1^T \Delta^k(x) dx-C_2^{(k)}T^{1+k/4}.$$

Lau and Tsang \cite{LT3} proved  $F_2(T)=\Omega(T\log^2 T).$ A weaker $\Omega$-result
of $F_2(T)$ can be found in \cite{IO}.
By using the sharpest known $\Omega$-result result of $\Delta(x)$  due to K. Soundararajan \cite{So}, Ivi\'c and Zhai \cite{IZ} proved that the estimate
$$F_k(T)=\Omega(T^{\frac{k+1}{4}}(\log T)^{\frac{k-3}{4}}
(\log\log T)^{\frac{3k+3}{4}(2^{4/3}-1)}
(\log\log\log T)^{-\frac{5k+5}{8}}) $$
holds for any $k\geq 3.$ However for odd $k\geq 3,$ the above $\Omega$-result of $F_k(T)$ can be improved to
$F_k(T)=\Omega(T^{\frac{k+2}{4}}(\log T)^{-5} )  $
(see Page 2567 of \cite{Ts2},  or \cite{TZ}).

{\bf Remark 1}. Zhai applied the above approach to study some similar error terms
in analytic number theory. Actually the approach can be applied to all error terms connected with the coefficients of Dirichlet series belonging to
the  Selberg class $\mathcal{S}$ of degree two. We mention the following examples.

(1) The power moments of $P(x)$ defined in (\ref{yi.1}) were investigated in \cite{Ka, Ts, Zw1,Zw2}.

(2) Let $a(n)$ denote the Fourier coefficients of a holomorphic cusp form of weight
$k= 2\ell \geq 12$ for the full modular group and    define
$A(x)=\sum_{n\leq x}a(n).$ The moments of $A(x)$ were investigated in \cite{Cai, Iv5, Zw1, Zw2}.

(3) Let $E(t)$ be defined by
$$E(t):=\int_0^t|\zeta(1/2+iu)|^2du-t\left(\log \frac{t}{2\pi}+2\gamma-1\right), \ \ t\geq 2.$$
The moments of $E(t)$ were investigated in \cite{BW, HI, HB, HB2, HBT,  Ju1, Ju2, Meu, Zw1, Zw2, Zw7}.

{\bf Remark 2.} The short interval results of moments of $\Delta(x)$ were studied in
 \cite{IS, LT, No,  Zw4}.

{\bf Remark 3.} Let $k\geq 1$ be a fixed integer. The discrete moments of the form $\sum_{n\leq x}\Delta^k(n) $ were studied by several authors.  See, for example,
\cite{CFTZ, Fu1, Fu2, Ha, Vo}.

{\bf Remark 4.}  Zhai's approach \cite{Zw2} was used by many authors to study moments
of error terms in analytic number theory. See, for example,
\cite{CLSXZ, FZ, HYZ, HTZZ, HWZ, LS, LS2, Lk, LW, SZZ, TZZ, TZ, Zw5, ZY}.

\subsection{Sign changes of $\Delta(x)$}\

By Cauchy's inequality we get
\begin{eqnarray*}
\left(\int_1^T\Delta^2(x)dx\right)^2&&\leq \int_1^T|\Delta(x)|dx\times \int_1^T|\Delta(x)|^3dx\\
&&=T^{7/4}\times \int_1^T|\Delta(x)|dx\times T^{-7/4}\int_1^T|\Delta(x)|^3dx
\end{eqnarray*}
which combining (\ref{er.4}) and (\ref{er.10}) with $A=3$ gives
$$\int_1^T|\Delta(x)|dx  \gg T^{5/4}.$$
Comparing this estimate and (\ref{er.3}), we see that $\Delta(x)$ changes sign  infinitely often when $x$ tends to infinity.

In 1955 Tong \cite{To1} proved the following fact:

{\it There exists two  positive  constants $c_1>0$ and $c_2>0$ such that,
for any $X\geq 2$ and  any $v\in [-c_1X^{1/4}, c_1X^{1/4} ],$ there is an
$x\in [X, X+c_2\sqrt X]$ for which $\Delta(x)=v.$ In particular,
there are two points $x_1, x_2\in [X, X+c_2\sqrt X]$ such that
$ \Delta(x_1)\geq c_3T^{1/4},\ \ \Delta(x_2)\leq -c_3T^{1/4}$ for some positive
constant $c_3>0.$ }

 Heath-Brown and Tsang \cite{HBT}   showed that the above length   $\sqrt X$ is almost optimal:

{\it Let $\delta$ be a sufficiently small positive number. Then for $X\geq X_0(\delta)$, there are at least $ \delta \sqrt{X}\log^{-5} X$ disjoint subintervals of length $\delta \sqrt{X}\log^{ 5} X$ in $[X, 2X]$ such that
$|\Delta(x)| > \delta x^{1/4} $ whenever $x$ lies in any of
these subintervals. In particular $\Delta(x)$ does not change sign in any of these subintervals.}

As was pointed out in \cite{Ts3}, $|\Delta(x)|$ in the above result can be replaced by $\Delta^{+}(x)=(|\Delta(x)|+\Delta(x))/2$ and $\Delta^{-}(x)=(|\Delta(x)|-\Delta(x))/2$ respectively (see also \cite{TZ}).

Ivi\'c and Zhai \cite{IZ4} generalized the above result of  Heath-Brown and Tsang to the short interval case. They proved that:

{\it Suppose $T,H $ are large parameters and $C > 1$ is a large constant
such that $CT^{3/4}\log\log T\leq  H \leq T$. Then in the interval $[T,T + H]$ there are $HT^{-1/2}\log^{5} T$  subintervals of length 	$T^{1/2}\log^{-5} T $ such that on each subinterval one has $\pm\Delta(x) \geq  c_{\pm}T^{1/4}$ for some $c_{\pm}> 0$.}

In \cite{Ts2} Tsang proposed the following

{\bf Problem 1.} Do there exist intervals $I = [X, X + X^\beta]$ with $\beta> 1/4$ such that the estimate
\begin{equation}\label{tsang-problem}
 \int_{I} |\Delta(x)|dx\ll T^{\beta+1/4-\delta}
 \end{equation}holds
for some small positive $\delta$?

Ivi\'c and Zhai \cite{IZ4} proved that Tsang's problem is not true for $1/2<\beta<1$. However, Tsang's problem for the remaining range $1/4<\beta\leq 1/2$ is very difficult, which is closely related to the number of sign changes of $\Delta(x)$. Ivi\'c and Zhai \cite{IZ4} proved that   if   $\Delta(x)$ could have enough sign changes, then (\ref{tsang-problem}) in this range would be true. Roughly speaking, they proved that  if $R\gg T^{1/2+\delta}, $ then   the estimate (\ref{tsang-problem}) holds, where $R$ denote the number of sign changes of
$\widetilde{\Delta}(x)$ in the interval $[T, 2T]$ and where
$$\widetilde{\Delta}(x)= \frac{x^{1/4}}{\sqrt 2 \pi} \sum_{n\leq T}\frac{d(n)}{n^{3/4}}
\cos(4\pi\sqrt {nx}- \pi/4).$$

\subsection{Hybrid moments of $\Delta(x)$ and $ \zeta(1/2 + ix)$}\

 Ivi\'c \cite{Iv4} proved several hybrid  results involving
the mean values of $\Delta(x)$, $E(t)$ and $\Delta^{*}(x),$ where
\begin{eqnarray*}
\Delta^{*}(x)&&:=-\Delta(x)=2\Delta(2x)-\frac 12\Delta(4x)\\
&& =\frac 12\sum_{n\leq 4x}(-1)^nd(n)-x(\log x+2\gamma-1)
\end{eqnarray*}
is the  "modified" divisor function, introduced and studied by M. Jutila
\cite{Ju1, Ju2}. For any integer $k\geq 1,$ define
$$\mathfrak{J}_k(T):=\int_1^T \left(\Delta^{*}\left(\frac{t}{2\pi}\right)\right)^k  |\zeta(1/2+it)|^{2}dt.$$
Ivi\'c  \cite{Iv4} proved the upper bound
\begin{eqnarray*}
&&\mathfrak{J}_{1}(T)\ll T^{7/6}\log^{7/2} T
\end{eqnarray*}
and proved the asymptotic formulas
\begin{eqnarray*}
&&\mathfrak{J}_{2}(T)=  T^{3/2}(c_1\log T+c_2)+O(T^{17/12+\varepsilon}),\\
&&\mathfrak{J}_{3}(T)=  T^{7/4}(c_3\log T+c_4)+O(T^{27/16+\varepsilon}),
\end{eqnarray*}
where $c_j(j=1, 2, 3, 4)$ are constants.

Ivi\'c and Zhai \cite{IZ2} studied the hybrid power moments of $\Delta(t)$ and
$|\zeta(1/2+it)|^2.$ For any fixed integers $k\geq 1$ and $m\geq 1, $ define
$$J_{k,m}(T):=\int_1^T \Delta^k(t)  |\zeta(1/2+it)|^{2m}dt.$$

 Ivi\'c and Zhai \cite{IZ2} proved
\begin{equation}\label{er.16}
J_{1,1}(T)\ll T(\log T)^4
\end{equation}
and they   conjectured that
\begin{equation}
 J_{1,1}(T)=\frac T4\left(\log \frac{T}{2\pi}+2\gamma-1\right)+O(T^{3/4+\varepsilon}).
\end{equation}

When $2\leq k\leq 8,$ Ivi\'c and Zhai \cite{IZ2} proved that the asymptotic formula
\begin{equation}\label{er.18}
J_{k,1}(T) =c_1(k)T^{1+\frac k4}\log T+c_2(k)T^{1+\frac k4}+O(T^{1+\frac k4-\eta_k+\varepsilon})
\end{equation}
holds, where $c_1(k)$ and $c_2(k)$ are explicit constants,
and $\eta_2=3/20, $   $\ \eta_3=\eta_4=1/10, $   $\ \eta_5=3/80, $  $\ \eta_6=35/4742, $   $\ \eta_7=17/6312, $ $\ \eta_8=8/9433. $

Ivi\'c   \cite{Iv3} studied the upper bound of $J_{k,m}(T) $ for $1\leq  k\leq 4$ and
$m=2,3.$
He proved  the following estimates
\begin{eqnarray*}
&&J_{1,2}(T)\ll T^{41/32}\log^C T,\ \ \ (41/32=1.28125),\\
&&J_{2,2}(T)\ll T^{25/16}\log^C T,\ \ \ (25/16=1.5625),\\
&&J_{3,2}(T)\ll T^{59/32}\log^C T,\ \ \ (59/32=1.84375),\\
&&J_{4,2}(T)\ll T^{17/8}\log^C T,\ \ \  \ (17/8=2.125),\\
&&J_{1,3}(T)\ll T^{19/32}\log^C T,\ \ \ (19/32=1.53125),\\
&&J_{2,3}(T)\ll T^{29/16}\log^C T,\ \ \ (29/16=1.8125),
\end{eqnarray*}
where $C>0$ is some fixed constant.
Ivi\'c   \cite{Iv3} even conjectured that
the asymptotic formula
\begin{eqnarray}\label{er.19}
&&J_{k,m}(T) =T^{1+\frac k4}Q_{m^2}(\log T)+O(T^{1+\frac k4-\rho_{k,m}+\varepsilon}),\ \ (k\geq 2)
\end{eqnarray}
holds, where $Q_{m^2}(x)$ is a polynomial of degree $m^2$, whose
coefficients depend on $k, m$ and where $\rho_{k,m}>0$ is a constant.

\section{Moments of   $\Delta_3(x)$ and some generalizations  }

\subsection{The general divisor problem}\

Let $k\geq 2$ be a fixed integer. The divisor function $d_k(n)$ counts the number of ways $n$ can be written as $k$ natural numbers and $d_2(n)=d(n)$. Define
$$D_k(x)=\sum_{n\leq x}d_k(n)=M_k(x)+\Delta_k(x),$$
where
$$M_k(x)=Res_{s=1}\ \zeta^k(s)x^ss^{-1}$$ is the main term
and $\Delta_k(x)$ is the error term. Studying properties of $\Delta_k(x)$
is an important problem in number theory(see, for example,  Ivi\'c \cite{Iv, Iv2, Iv6}).

The upper bound of $\Delta_3(x)$ was studied by many authors and the latest result
of Kolesnik \cite{kol} reads
\begin{equation}
\Delta_3(x)\ll x^{43/96+\varepsilon}.
\end{equation}

Tong \cite{To2} proved the mean square result
\begin{equation}\label{san.2}
\int_1^T\Delta_3^2(x)dx=C_3^{(2)}T^{5/3}+O(T^{14/9+\varepsilon}),
\end{equation}
where $C_3^{(2)}>0$ is a positive constant. His main ingredient is to replace $\Delta_3(x)$ by a certain integral such that the difference between $\Delta_3(x)$ and this
integral is very small on average.

Heath-Brown \cite{HB} proved that
\begin{equation}
\int_1^T|\Delta_3(x)|^3 dx\ll  T^{2+\varepsilon}.
\end{equation}

Heath-Brown \cite{HB} also  proved that for any fixed real number $0<k<3$ there holds
\begin{equation}\label{san.4}
\int_1^T|\Delta_3(x)|^k dx=T^{1+\frac k3}(c_k+o(1))
\end{equation}
for some constant $c_k>0.$

When $k\geq 4,$ the upper bound and mean square results of $\Delta_k(x)$,  can see, for example, Heath-Brown \cite{HB3}, Ivi\'c \cite{Iv, Iv2} and Zhang \cite{Zh}.

\subsection{ Sign changes of $\Delta_3(x)$ }\

Suppose $k\geq 3$ is a fixed integer. Tong \cite{To1} proved that:
 {\it There exist two positive constants $c_k$ and $C_k$, such that for every $T\geq 2$ and every $v$ lying in the interval $[-c_kT^{(k-1)/2k}, c_kT^{(k-1)/2k}],$  the equation $\Delta_k(y)=v$ has  at least  one solution lying in the interval
$[T, T+C_k T^{1-1/k}]$.}

The above result shows that $\Delta_k(x)$ must have at least one sign change in the interval $[T, T+C_k T^{1-1/k}] $ for any $T\geq 2.$  Especially $\Delta_3(x)$ must have at least one sign change  in the interval $[T, T+C_3T^{2/3}]. $ Note that
here $2/3$ is   best possible(see Remark 18 in \cite{CTZ2}).

In the opposite direction, Cao, Tanigawa and Zhai \cite{CTZ2} proved the following result:
{\it Let $T$ be a sufficiently large real number. Then in the interval $[T, 2T]$, there are at least $ T^{1/2-\varepsilon}$ disjoint subintervals of length
$T^{1/2-\varepsilon}$   such that
$|\Delta_3(x)| > c T^{1/3} $ whenever $x$ lies in any of
these subintervals, where $c>0$ is a positive number. In particular $\Delta_3(x)$ does not change sign in any of these subintervals.}

\subsection{ A generalization of Tong's approach}\

Let $\mathcal{S}_{real}$ denote a subset of the Selberg class $\mathcal{S}$ such that
$\mathcal{S}_{real}$ includes all Dirichlet series $\mathfrak{L}_f(s)$ defined by (\ref{yi.3}) with $f(n)\in \mathbb{R} \ (n\geq 1).$

Suppose $\{f(n)\}(n=1,2,3,\cdots)$ is a sequence of real numbers such that
its corresponding Dirichlet series $ \mathfrak{L}_f(s)\in {\mathcal S}_{real} $ is a function of degree $d.$ Let $1/2\le \sigma^{*}<1 $ denote the infimum
of $\sigma$ such that
$$\int_0^T| \mathfrak{L}_f(\sigma+it)|^2dt\ll T^{1+\varepsilon}$$
 holds for any $\varepsilon>0.$ We suppose that $\sigma^{*}$ satisfies the condition
\begin{equation}\label{sigma}
\sigma^{*}<(d+1)/2d,
\end{equation}

Cao, Tanigawa and Zhai \cite{CTZ3} generalized Tong's result to the general situation.
They proved that: {\it
 Suppose that $d\ge 5/3$ and (\ref{sigma}) holds, then
we have
\begin{equation}
\int_1^T E_f^2(x)dx=C_fT^{2-1/d}+O(T^{2-\frac{3-4\sigma^{*}}{2d(1-\sigma^{*})-1}+\varepsilon}),
\end{equation}
where $C_f>0$ is a positive constant.}

As a corollary, they deduced the following results.

(1) {\it If $\mathfrak{L}_f(s)\in \mathcal{S}_{real}$ is a function of degree $2$, then
\begin{equation*}
\int_1^TE_f^2(x)dx=C_fT^{3/2}+O(T^{1+\varepsilon}).
\end{equation*}}

(2){\it  If $\mathfrak{L}_f(s)\in \mathcal{S}_{real}$ is a function of degree $3$
such that $\mathfrak{L}_f(s)=\mathfrak{L}_1(s)\mathfrak{L}_2(s)$, where
$ \mathfrak{L}_1(s) \in \mathcal{S}_{real}$ is a function of degree $1$
and $ \mathfrak{L}_2(s) \in \mathcal{S}_{real}$ is a function of degree $2,$
then
\begin{equation*}
\int_1^TE_f^2(x)dx=C_fT^{5/3}+O(T^{8/5+\varepsilon}).
\end{equation*}
Furthermore if we assume that
\begin{equation*}
\int_1^T  |\mathfrak{L}_2(1/2+it)|^6\ll T^{2+\varepsilon},
\end{equation*}
then
\begin{equation*}
\int_1^TE_f^2(x)dx=C_fT^{5/3}+O(T^{14/9+\varepsilon}).
\end{equation*}
}

\section{Hybrid moments of   remainders in   divisor problems}

\subsection{Some hybrid Moments of   $\Delta_2(x), \Delta_3(x)$ and $\Delta_4(x)$  }\

  Ivi\'c and   Zhai \cite{IZ3} first studied the hybrid Moments of
 $\Delta_k(x)\ (k\geq 2).$ They proved that the estimates
\begin{eqnarray}\label{si.1}
&&\int_1^T\Delta_2(x)\Delta_3(x)dx\ll T^{13/9}\log^{10/3} T
\end{eqnarray}
and
\begin{eqnarray*}
&&\int_1^T\Delta_2(x)\Delta_4(x)dx\ll T^{25/16+\varepsilon}
\end{eqnarray*}
hold.

Zhai \cite{Zw6} proved that
\begin{eqnarray*}
&&\int_1^T\Delta_2^2(x)\Delta_3(x) dx\ll T^{16/9+\varepsilon},\ \
\int_1^T\Delta_2^3(x)\Delta_3(x)dx\ll T^{73/36+\varepsilon}.
\end{eqnarray*}

Cao, Tanigawa and Zhai \cite{CTZ} proved that the estimate
\begin{eqnarray*}
&&\int_1^T\Delta_3^2(x)\Delta_2(x) dx\ll T^{67/36+\varepsilon}
\end{eqnarray*}
holds. They also proved that the asymptotic formulas
\begin{eqnarray}\label{si.2}
&&\int_1^T\Delta_3^2(x)\Delta_2^2(x) dx=\frac{15}{13}C_2^{(2)}C_3^{(2)}T^{13/6}+O( T^{13/6-2/45+\varepsilon})
\end{eqnarray}
and
\begin{eqnarray}\label{si.3}
&&\int_1^T\Delta_3^2(x)\Delta_2^3(x) dx=\frac{70}{29\sqrt 2}C_2^{(3)}C_3^{(2)}T^{29/12}+O( T^{29/12- 178/16647+\varepsilon})
\end{eqnarray}
hold, where $C_2^{(2)},   C_2^{(3)}, C_3^{(2)}$ were defined in (\ref{er.4}), (\ref{er.8}), (\ref{san.2}) respectively.

Recently, Cai and Zhai \cite{CZ} proved
\begin{eqnarray}\label{si.4}
&&\int_1^T\Delta_2(x)\Delta_3(x)\Delta_4(x)dx\ll T^{ 93/48+\varepsilon},\ \ \
\int_1^T\Delta_3(x)\Delta_4(x)dx\ll T^{5/3+\varepsilon}.
\end{eqnarray}

\subsection{Joint distribution of signs of $\Delta_2(x) $ and $\Delta_3(x)$ }\

In this subsection we give some results about the joint distribution of signs of $\Delta_2(x) $ and $\Delta_3(x)$.

{\bf Theorem 1.} {\it Suppose $T$ is a large real number. Then $\Delta_2(x)\Delta_3(x)$ changes signs in the interval $[T, T+T^{139/144+\varepsilon}],$ where $  139/144=0.9652\overline{7}.$ }

\begin{proof}
Let $ T^{139/144+\varepsilon}\leq H\leq T.$ By (\ref{si.1}) we have
\begin{eqnarray}\label{si.5}
&&\int_T^{T+H}\Delta_2(x)\Delta_3(x)dx\ll T^{13/9}\log^{10/3} T.
\end{eqnarray}

From (\ref{si.2}) we get
\begin{eqnarray}\label{si.6}
\ \ \ \ \int_T^{T+H}\Delta_3^2(x)\Delta_2^2(x) dx
=C_0((T+H)^{\frac{13}{6}}-T^{\frac{13}{6}})+O( T^{\frac{13}{6}- \frac{2}{45}+\varepsilon/2})  \gg HT^{7/6},
\end{eqnarray}
where
$C_0=\frac{15}{13}C_2^{(2)}C_3^{(2)}.$

Suppose $6<k_2\leq 9 $ and $2<k_3<3$ are fixed real numbers such that
$1-2/k_2-2/k_3>0.$ Let
$a: =\frac{1-2/k_2-2/k_3}{1-1/k_2-1/k_3},$ we have
$$a+\frac{2-a}{k_2}+\frac{2-a}{k_3}=1,\ \ \ 2-a=\frac{1}{1-1/k_2-1/k_3}.$$
By H\"{o}lder's inequality  we have
\begin{eqnarray*}
&&\ \ \ \int_T^{T+H}\Delta_3^2(x)\Delta_2^2(x) dx   =
\int_T^{T+H}|\Delta_3(x)\Delta_2(x)|^{a}|\Delta_3(x)|^{2-a}|\Delta_2(x)|^{2-a} dx \\
&&\leq  \left(\int_T^{T+H}|\Delta_3(x)\Delta_2(x)|   dx\right)^{a}  \times
\left( \int_T^{T+H}  |\Delta_2(x)|^{k_2} dx  \right)^{\frac{2-a}{k_2}}\nonumber\\
&&\ \ \ \  \times
\left(\int_T^{T+H} |\Delta_3(x)|^{k_3}  dx\right)^{\frac{2-a}{k_3}} \nonumber\\
&&=  \left(\int_T^{T+H}|\Delta_3(x)\Delta_2(x)|   dx\right)^{\frac{1-2/k_2-2/k_3}{1-1/k_2-1/k_3}}  \times
\left( \int_T^{T+H}  |\Delta_2(x)|^{k_2} dx  \right)^{\frac{1/k_2}{1-1/k_2-1/k_3 }}\nonumber\\
&&\ \ \ \  \times
\left(\int_T^{T+H} |\Delta_3(x)|^{k_3}  dx\right)^{\frac{1/k_3}{1-1/k_2-1/k_3}} \nonumber.
\end{eqnarray*}
Thus we get
\begin{eqnarray}\label{si.7}
&&\int_T^{T+H}|\Delta_3(x)\Delta_2(x)|   dx \geq
   \left(\int_T^{T+H}\Delta_3^2(x)\Delta_2^2(x) dx\right)^{\frac{1-1/k_2-1/k_3}{1-2/k_2-2/k_3} }\\
&& \ \ \ \  \times
\left( \int_T^{T+H}  |\Delta_2(x)|^{k_2} dx  \right)^{-\frac{1/k_2}{1-2/k_2-2/k_3 }}   \times
\left(\int_T^{T+H} |\Delta_3(x)|^{k_3}  dx\right)^{-\frac{1/k_3}{1-2/k_2-2/k_3}} \nonumber.
\end{eqnarray}

From (\ref{san.4}) we get
\begin{eqnarray}\label{si.8}
&&  \int_T^{T+H} |\Delta_3(x)|^{k_3}  dx\leq
\int_T^{2T} |\Delta_3(x)|^{k_3}  dx\ll T^{1+k_3/3}=\frac TH\times HT^{k_3/3}.
\end{eqnarray}
From results in \cite{Zw4}   we get
\begin{eqnarray}\label{si.9}
&&  \int_T^{T+H} |\Delta_2(x)|^{k_2}  dx \asymp  HT^{\frac{k_2}{4}}.
\end{eqnarray}

From (\ref{si.6})-(\ref{si.9}) we get
\begin{eqnarray}\label{si.10}
&&\int_T^{T+H}|\Delta_2(x)\Delta_3(x)|dx\gg H T^{7/12}\times
(TH^{-1})^{-\frac{1/k_3}{1-2/k_2-2/k_3}}.
\end{eqnarray}
Take $k_2=9$  and $1/k_3=1/3+\delta$,  where $\delta>0$ is a very small number
such that
$$3<\frac{1/k_3}{1-2/k_2-2/k_3}<3+\varepsilon. $$
From (\ref{si.10}) we get
\begin{eqnarray}\label{si.11}
&&\int_T^{T+H}|\Delta_2(x)\Delta_3(x)|dx\gg H T^{7/12}\times
(TH^{-1})^{-( 3+\varepsilon)}\gg H^4T^{-29/12-\varepsilon}.
\end{eqnarray}
Now Theorem 1 follows from (\ref{si.5}) and (\ref{si.11}).
\end{proof}

{\bf Lemma 4.1.} {\it Suppose $2\leq U\leq \sqrt{T}/2$, $1/4<\theta_2<1/3$ and $1/3<\theta_3<1/2$ are real numbers such that
$\Delta_2(x)\ll x^{\theta_2},\ \ \Delta_3(x)\ll x^{\theta_3}.$ Then we have
\begin{eqnarray*}
&&\ \ \ \int_T^{2T}
\max_{0\leq u\leq U}(\Delta_2(x+u)\Delta_3(x+u)-\Delta_2(x)\Delta_3(x) )^2dx\\&&
\ll T^{1+2\theta_3}U\log^5 T  + T^{4/3+2\theta_2+\varepsilon}U^{1/2}+T^{19/15+2\theta_2+\varepsilon}U^{4/5} .
\end{eqnarray*}
}

\begin{proof}
For any $T\leq x\leq 2T,$ we have
\begin{eqnarray*}
&&\ \ \ \ \ \ \Delta_2(x+u)\Delta_3(x+u)-\Delta_2(x)\Delta_3(x)\\
&&=(\Delta_2(x+u) -\Delta_2(x)) \Delta_3(x+u)
 +\Delta_2(x)(\Delta_3(x+u)- \Delta_3(x))\\
 &&\ll T^{\theta_3} \times |\Delta_2(x+u) -\Delta_2(x)|
 +T^{\theta_2}\times |\Delta_3(x+u)- \Delta_3(x)|,
\end{eqnarray*}
which implies that
\begin{eqnarray*}
&&\ \ \ \ \ \ \max_{0\leq u\leq U}(\Delta_2(x+u)\Delta_3(x+u)-\Delta_2(x)\Delta_3(x))^2\\
&&\ll  T^{2\theta_3} \times \max_{0\leq u\leq U}(\Delta_2(x+u) -\Delta_2(x))^2
 +T^{2\theta_2}\times \max_{0\leq u\leq U}(\Delta_3(x+u)- \Delta_3(x))^2.
\end{eqnarray*}
So we have
 \begin{eqnarray}\label{si.12}
&&\ \ \ \ \ \ \ \ \ \ \ \  \int_T^{2T}\max_{0\leq u\leq U}(\Delta_2(x+u)\Delta_3(x+u)-\Delta_2(x)\Delta_3(x))^2dx\\
&&\ll  T^{2\theta_3}   \int_T^{2T}\max_{0\leq u\leq U}(\Delta_2(x+u) -\Delta_2(x))^2dx
 +T^{2\theta_2} \int_T^{2T}\max_{0\leq u\leq U}(\Delta_3(x+u)- \Delta_3(x))^2dx.\nonumber
\end{eqnarray}

Theorem 4 of \cite{IZ4} reads
\begin{eqnarray}\label{si.13}
\int_T^{2T}\max_{0\leq u\leq U}(\Delta_2(x+u) -\Delta_2(x))^2dx\ll TU\mathcal{L}^5.
\end{eqnarray}

Since $U\ll T^{1/2},$   Lemma 19 of \cite{CTZ2} reads
\begin{eqnarray}\label{si.14}
\int_T^{2T}\max_{0\leq u\leq U}(\Delta_3(x+u) -\Delta_3(x))^2dx\ll T^{4/3+\varepsilon}U^{1/2}+T^{19/15+\varepsilon}U^{4/5} .
\end{eqnarray}

Now Lemma 4.1 follows from (\ref{si.12})-(\ref{si.14}).
\end{proof}

{\bf Theorem 2.} {\it Suppose $T$ is a large real number and $\delta>0$ is a fixed small number.

(1) There are  $ \gg T^{35/48+\varepsilon}$ disjoint subintervals of length $\gg T^{13/48-\varepsilon}$ in $[T, 2T]$ such that whenever $x$ lies in any of these subintervals, we have
\begin{eqnarray}\label{si.15}
| \Delta_2(x)|>\delta x^{1/4}, \ \ |\Delta_3(x)| > \delta x^{1/3}, \ \ sign(\Delta_2(x)\Delta_3(x))=\pm 1.
\end{eqnarray}

(2) If $\Delta_2(x)\ll x^{1/4+\varepsilon}, \ \Delta_3(x)\ll x^{1/3+\varepsilon}, $ then there are  $ \gg T^{1/2+\varepsilon}$ disjoint subintervals of length $\gg T^{1/2-\varepsilon}$ in $[T, 2T]$ such that (\ref{si.15}) holds whenever $x$ lies in any of these subintervals. }

\begin{proof}
{\bf Step 1.} From (\ref{si.11}) with $H=T$ we get
\begin{eqnarray}\label{si.16}
&&\int_T^{2T}|\Delta_2(x)\Delta_3(x)|dx\gg  T^{19/12}.
\end{eqnarray}
From (\ref{si.1})  we get
\begin{eqnarray}\label{si.17}
&&\int_T^{2T} \Delta_2(x)\Delta_3(x) dx\ll  T^{13/9}\log^{10/3} T.
\end{eqnarray}

{\bf Step 2.} Define
\begin{eqnarray}\label{si.18}
\mathfrak{E}_{+}(x)=\left\{\begin{array}{ll}
\Delta_2(x)\Delta_3(x) ,&\mbox{if $\Delta_2(x)\Delta_3(x)>0,$}\\
0,&\mbox{if $\Delta_2(x)\Delta_3(x)\leq 0,$}
\end{array}\right.
\end{eqnarray}
and
\begin{eqnarray}\label{si.19}
\mathfrak{E}_{-}(x)=\left\{\begin{array}{ll}
-\Delta_2(x)\Delta_3(x) ,&\mbox{if $\Delta_2(x)\Delta_3(x)<0,$}\\
0,&\mbox{if $\Delta_2(x)\Delta_3(x)\geq 0.$}
\end{array}\right.
\end{eqnarray}

It is easily seen that
\begin{eqnarray*}
\int_T^{2T}|\Delta_2(x)\Delta_3(x)|dx=2\int_T^{2T}\mathfrak{E}_{\pm}(x)dx
\mp\int_T^{2T} \Delta_2(x)\Delta_3(x) dx.
\end{eqnarray*}

Thus from (\ref{si.2}), (\ref{si.18}), (\ref{si.19}) and Cauchy's inequality we get
\begin{eqnarray*}
 T^{19/12 }\ll \int_T^{2T}|\Delta_2(x)\Delta_3(x)|dx
\ll T^{1/2}\left(\int_T^{2T}\mathfrak{E}_{\pm}^2(x)dx\right)^{1/2}.
\end{eqnarray*}
So we get that the estimate
\begin{eqnarray}\label{si.20}
\int_T^{2T}\mathfrak{E}_{\pm}^2(x)dx\geq    c_{\pm}T^{13/6}
\end{eqnarray}
holds for some positive constant $c_{\pm}>0.$

{\bf Step 3.}  Let $\delta>0$ be a small constant. Suppose
\begin{eqnarray*}
 G(x)&&=(|\Delta_2(x)|-\delta x^{1/4 })\times
(|\Delta_3(x)|-\delta x^{1/3 })\\
&&=|\Delta_2(x) \Delta_3(x)|-\delta x^{1/3 }|\Delta_2(x)|
-\delta x^{1/4 }|\Delta_3(x)|+\delta^2 x^{7/12 }.
\end{eqnarray*}

Define
\begin{eqnarray*}
\mathfrak{G}_{+}(x)=\left\{\begin{array}{ll}
 G(x),&\mbox{if $\Delta_2(x)\Delta_3(x)>0,$}\\
0,&\mbox{if $\Delta_2(x)\Delta_3(x)\leq 0,$}
\end{array}\right.
\end{eqnarray*}
and
\begin{eqnarray*}
\mathfrak{G}_{-}(x)=\left\{\begin{array}{ll}
G(x) ,&\mbox{if $\Delta_2(x)\Delta_3(x)<0,$}\\
0,&\mbox{if $\Delta_2(x)\Delta_3(x)\geq 0.$}
\end{array}\right.
\end{eqnarray*}

Let $I=[T, 2T]$ and write $I=I_{+}\bigcup I_{-}\bigcup I_{0},$ where
\begin{eqnarray*}
&&I_{+}=\{x\in I:  \Delta_2(x)\Delta_3(x)>0\},\\
&&I_{-}=\{x\in I:  \Delta_2(x)\Delta_3(x)<0\},\\
&&I_{0}=\{x\in I:  \Delta_2(x)\Delta_3(x)=0\}.
\end{eqnarray*}
We will show that  if $\delta>0$ small enough then
\begin{eqnarray}\label{si.21}
\int_{T}^{2T}\mathfrak{G}_{\pm}^2(x)dx\geq   \frac{c_{\pm} }{2} T^{13/6 }.
\end{eqnarray}

Write
\begin{eqnarray}\label{4.22}
G^2(x)=|\Delta_2(x)\Delta_3(x)|^2+\sum_{j=1}^{8}G_j(x),
\end{eqnarray}
where
\begin{eqnarray*}
&&G_1(x)=\delta^2 x^{2/3 } \Delta_2^2(x), \ \ \
 G_2(x)=\delta^2 x^{1/2 }\Delta_3^2(x),\ \ \
 G_3(x)=\delta^4 x^{7/6},\\
&&G_4(x)=-2\delta x^{1/3 }|\Delta_2^2(x) \Delta_3(x)| ,\ \
G_5(x)= -2  \delta x^{1/4 }  |\Delta_2(x) \Delta_3^2(x)|,\\
&&G_6(x)=4\delta^2 x^{7/12 }|\Delta_2(x) \Delta_3(x)|,\ \
G_7(x)=-2\delta^3 x^{11/12}|\Delta_2(x)|,\\
&&
G_8(x)=-2\delta^3 x^{10/12 }|\Delta_3(x)|.
\end{eqnarray*}

With the help of moments results of $\Delta_2(x)$ and $\Delta_3(x)$ we get
\begin{eqnarray*}
\int_T^{2T}G_j(x)dx\ll  \delta T^{13/6}, \ \ \ 1\leq j\leq 8.
\end{eqnarray*}

So there is a constant $C=C(\varepsilon)>0$ such that
\begin{eqnarray}\label{si.22}
 \sum_{j=1}^8 \left|\int_T^{2T}G_j(x)dx\right| \leq \delta C T^{13/6}
 \leq \frac{c_{\pm}}{2}T^{13/6}
\end{eqnarray}
if $0<\delta<c_{\pm}/2C.$

From (\ref{4.22}) and (\ref{si.22}) we get
\begin{eqnarray}\label{si.23}
\ \  \ \ \ \int_{T}^{2T}\mathfrak{G}_{\pm}^2(x)dx&&=\int_{I_{\pm}}G^2(x)dx
=\int_{I_{\pm}} |\Delta_2(x)\Delta_3(x)|^2dx+\sum_{j=1}^{8}\int_{I_{\pm}} G_j(x)dx \\
&&=\int_{T}^{2T}\mathfrak{E}_{\pm}^2(x)dx+\sum_{j=1}^{8}\int_{I_{\pm}} G_j(x)dx
\geq \frac{ c_{\pm}}{2}T^{13/6}.\nonumber
\end{eqnarray}
Namely, (\ref{si.21}) holds.

{\bf Step 4.} Suppose $T^\varepsilon\ll U\ll T^{1-\varepsilon}.$
Define
\begin{eqnarray*}
\omega_{\pm}(t)=\mathfrak{G}_{\pm}^2(t)-4\max_{0\leq u\leq U}(G(t+u)-G(t))^2-
\frac{c_{\pm}}{100}t^{7/6}.
\end{eqnarray*}

If $\omega_{\pm}(t)>0,$ then it follows that
\begin{eqnarray}\label{si.24}
&&|G(t)|\geq 2\max_{0\leq u\leq U}|G(t+u)-G(t)|,\ \ \ t\in I_{\pm},
\end{eqnarray}
\begin{eqnarray}\label{si.25}
&&  |G(t)|\geq    \frac{\sqrt{c_{\pm}}}{10}t^{7/12}.
\end{eqnarray}

From (\ref{si.24}) we can see that $G(x)$ does not change its sign on the
interval $[t, t + U]$, and if $G(x) \leq  0$ in this interval it contradicts with (\ref{si.25}).
Hence if $\omega_{\pm}(t)> 0$ we can conclude that
$G(x)>0 $ for any $x \in [t, t + U]$, which implies that either (\ref{si.15}) holds
  or
$$ |\Delta_2(x)|<\delta x^{1/4 }, \
 |\Delta_3(x)|<\delta x^{1/3 }, \ sign(\Delta_2(x)\Delta_3(x))=\pm 1. $$
If the latter case holds, then
$$ -\delta x^{1/4 }<|\Delta_2(x)|-\delta x^{1/4 }<0,\
-\delta x^{1/3 } <|\Delta_3(x)|-\delta x^{1/3 }<0, $$
which implies that $G(x)<\delta^2 x^{7/12}.$  From (\ref{si.25}) we see this is a contradiction  if
$\delta  $ is small enough.

{\bf Step 5.} So we are going to study the set
$$  \mathfrak{X} =\{t\in [T, 2T]: \ \omega_{\pm}(t)>0 \}. $$
From (\ref{er.10}) with $A=9,$ (\ref{san.4}) with $k=2.7$ and H\"{o}lder's inequality we get
\begin{eqnarray}\label{si.26}
&&\ \ \ \ \ \int_T^{2T}\omega_{\pm}(x)dx\leq \int_{\mathfrak{X}  }\omega_{\pm}(x)dx\leq
\int_{\mathfrak{X}  }  \Delta_2^2(x) \Delta_3^2(x) dx\\
&&\ll |\mathfrak{X}|^{\frac{1}{27}}\left(\int_{T}^{2T}|\Delta_2(x)|^{9}dx\right)^{\frac{2}{9}}
\left(\int_{T}^{2T}|\Delta_3(x)|^{2.7}dx\right)^{\frac{20}{27}}\nonumber\\
&&\ll |\mathfrak{X}|^{\frac{1}{27}}T^{\frac{115}{54}}\nonumber.
\end{eqnarray}

Now we prove a lower bound of $\int_T^{2T}\omega_{\pm}(x)dx.$ Recall
$$G(x)=|\Delta_2(x) \Delta_3(x)|-\delta x^{1/3 }|\Delta_2(x)|
-\delta x^{1/4 }|\Delta_3(x)|+\delta^2 x^{7/12 },$$
we have
\begin{eqnarray*}\ \ \ \ \ \
G(x+u)&&=|\Delta_2(x+u) \Delta_3(x+u)|-\delta (x+u)^{1/3 }|\Delta_2(x+u)|\\
&&\ \ \ \ \ -\delta (x+u)^{1/4 }|\Delta_3(x+u)|+\delta^2 (x+u)^{7/12 }\nonumber\\
&&=|\Delta_2(x+u) \Delta_3(x+u)|-\delta x^{1/3 }|\Delta_2(x+u)|\nonumber\\
&&\ \ \ \ \ -\delta x^{1/4 }|\Delta_3(x+u)|+\delta^2 x^{7/12 }+ \nonumber\\&&
\ \ \ \ +O(\delta ux^{-2/3 }|\Delta_2(x+u)|+\delta ux^{-3/4 }|\Delta_3(x+u)|+
\delta^2 ux^{-5/12 }).\nonumber
\end{eqnarray*}
So we have (suppose $U<T^{1/2}$)
\begin{eqnarray*}
&&\ \ \ \ \max_{0\leq u\leq U}(G(x+u)-G(x))^2\\
&&\ll \max_{0\leq u\leq U}(\Delta_2(x+u)\Delta_3(x+u)-\Delta_2(x)\Delta_3(x))^2\\
&&+\delta^2 x^{2/3 }\max_{0\leq u\leq U}(|\Delta_2(x+u)| -|\Delta_2(x)|)^2
+\delta^2 x^{1/2 }\max_{0\leq u\leq U}(|\Delta_3(x+u)| -|\Delta_3(x)|)^2\\
&&+\delta^2 U^2x^{-4/3 }\max_{0\leq u\leq U} \Delta_2^2(x+u) +\delta^2 U^2x^{-3/2 } \max_{0\leq u\leq U}\Delta_3^2(x+u) +
\delta^4 U^2x^{-5/6 }\\
&&\ll \max_{0\leq u\leq U}(\Delta_2(x+u)\Delta_3(x+u)-\Delta_2(x)\Delta_3(x))^2\\
&&+(\delta^2 x^{2/3 }+\delta^2 U^2x^{-4/3 })\max_{0\leq u\leq U}(|\Delta_2(x+u)| -|\Delta_2(x)|)^2\\
&& +(\delta^2 x^{1/2 }+\delta^2 U^2x^{-3/2 })\max_{0\leq u\leq U}(|\Delta_3(x+u)| -|\Delta_3(x)|)^2\\
&&+\delta^2 U^2x^{-4/3 } \Delta_2^2(x) +\delta^2 U^2x^{-3/2 }  \Delta_3^2(x) +
\delta^4 U^2x^{-5/6 }\\
&&\ll \max_{0\leq u\leq U}(\Delta_2(x+u)\Delta_3(x+u)-\Delta_2(x)\Delta_3(x))^2\\
&&+ \delta^2 x^{2/3 } \max_{0\leq u\leq U}(|\Delta_2(x+u)| -|\Delta_2(x)|)^2
  + \delta^2 x^{1/2 } \max_{0\leq u\leq U}(|\Delta_3(x+u)| -|\Delta_3(x)|)^2\\
&&+\delta^2 U^2x^{-4/3 } \Delta_2^2(x) +\delta^2 U^2x^{-3/2 }  \Delta_3^2(x) +
\delta^4 U^2x^{-5/6}.
\end{eqnarray*}

Thus from Lemma 4.1 with $\theta_2=131/416+\varepsilon,\ \theta_3=43/96+\varepsilon/4$, (\ref{si.13}) and (\ref{si.14}) we get that
\begin{eqnarray}\label{si.27}
&& \ \ \ \ \ \int_T^{2T}\max_{0\leq u\leq U}(G(x+u)-G(x))^2dx \\
&&\ll T^{1+2\theta_3}U\log^5 T  + T^{4/3+2\theta_2+\varepsilon}U^{1/2}+T^{19/15+2\theta_2+\varepsilon}U^{4/5}\ll T^{13/6-\varepsilon/2}  \nonumber
\end{eqnarray}
for $U: =T^{13/48-\varepsilon}.$

From Lemma 4.1 with $\theta_2=1/4+\varepsilon,\ \theta_3=1/3+\varepsilon/4$, (\ref{si.13}) and (\ref{si.14}) we obtain
\begin{eqnarray}\label{sss}
\int_T^{2T}\max_{0\leq u\leq U}(G(x+u)-G(x))^2dx \ll T^{13/6-\varepsilon/2}, \ \ \
 U=T^{1/2-\varepsilon}.
\end{eqnarray}

From (\ref{si.21}) and (\ref{si.27}) we get
\begin{eqnarray}\label{si.28}
&&\ \ \ \ \ \int_T^{2T}\omega_{\pm}(x)dx\\
&&=\int_T^{2T}\mathfrak{G}_{\pm}^2(x)dx-
4\int_T^{2T}\max_{0\leq u\leq U}(G(x+u)-G(x))^2dx-
\int_T^{2T}\frac{c_{\pm}}{100}x^{7/6}dx\nonumber\\
&&\geq \left(\frac{c_{\pm}}{2}-\frac{c_{\pm}}{100}\times \frac{6(2^{13/6}-1)}{13}\right)T^{13/6}+O(T^{13/6-\varepsilon})\nonumber\\
&&\geq  \frac{c_{\pm}}{4} T^{13/6} \nonumber.
\end{eqnarray}

From (\ref{si.26}) and (\ref{si.28}) we get
\begin{eqnarray}\label{si.29}
 |\mathfrak{X}|\gg T.
\end{eqnarray}
By (\ref{si.25}), (\ref{si.29}) and the value  $U=T^{13/48-\varepsilon}$ we get the first assertion of Theorem 2. Using (\ref{sss}) instead of (\ref{si.27}) in the above argument we get the second assertion of Theorem 2.
\end{proof}

\section{Some conjectures  }

In this section we propose some conjectures. Some special cases of these conjectures are known results. We recall that $\mathcal{S}$ denotes the Selberg class,   $\mathcal{S}_{real}$ denotes
a subset of $\mathcal{S}$ which includes all Dirichlet series (\ref{yi.3}) such that
$f: \mathbb{N}\rightarrow \mathbb{R}.$

\subsection{Conjectures on moments of one error term }\

{\bf Conjecture 1.} {\it Let $\ell\geq 2$ be a fixed integer. If $\mathfrak{L}_f(s)\in \mathcal{S}$ is a function of degree $d\geq 2$, then there exist constants
$\delta_1>0$ and $\delta_2>0$ such that the estimate
\begin{equation}\label{C-1}
\int_1^TE_f(x)dx\ll  T^{1+\frac{ (d-1)}{2d}-\delta_1+\varepsilon}
\end{equation}
and the asymptotic formula
\begin{equation}\label{C-2}
\int_1^TE_f^\ell(x)dx=C_{f,d}^{(\ell)}T^{1+\frac{\ell (d-1)}{2d}}+O(T^{1+\frac{\ell (d-1)}{2d}-\delta_2+\varepsilon})
\end{equation}
hold, where   $C_{f,d}^{(\ell)}\not= 0$ is a constant.}

{\bf Remark 1.} The estimate (\ref{C-1}) can be investigated by moment results of $\mathfrak{L}_f(s)$ and even sharper results are available when $d$ is small.  When $d=2,$ by Zhai's approach \cite{Zw2} we see that (\ref{C-2}) holds for $2\leq \ell\leq 7$ for any $\mathfrak{L}_f(s)\in \mathcal{S}.$ Especially when $\mathfrak{L}_f(s)$ can be written as a product of two functions of degree $1,$ then (\ref{C-2}) holds for $2\leq \ell\leq 9.$ When $d=3,$ in \cite{CTZ3, CTZ4} there are some examples for which  (\ref{C-2}) holds. When $d\geq 4,$ no examples are known for which  (\ref{C-2}) holds.

From Theorem 6 of \cite{CTZ4}, we have the following result: {\it If  $\mathfrak{L}_f(s)\in \mathcal{S}_{real}$ is a function of degree $d\geq 2$ and    $T\geq 10$ is a large real number. Then there exist $x_j\in [T, T+C(f) T^{1-1/d}](j=1,2)$ such that $$E_f(x_1)>c(f)T^{(d-1)/2d}, \ \ E_f(x_2)<-c(f)T^{(d-1)/2d},$$  where $C(f)$ and $c(f)$ are two positive constants.
Especially   $E_f(x)$ must have at least one sign change in the interval $[T, T+C(f) T^{1-1/d}] $. }

We remark that the exponent $1-1/d$ mentioned above {\it should be}   best possible. We have the following

{\bf Conjecture 2.} {\it   Let $T$ be a sufficiently large real number. Then in the interval $[T, 2T]$, there are at least $ T^{1/d+\varepsilon}$ disjoint subintervals of length
$T^{1-1/d-\varepsilon}$   such that
$\pm E_f(x) > c(f) T^{(d-1)/2d} $ whenever $x$ lies in any of
these subintervals, where $c(f)>0$ is a positive number.  }

{\bf Remark 2.} By the approach of \cite{HBT} we see that Conjecture 2 is true for $d=2.$ When the Lindel\"{o}f hypothesis of $\zeta(s)$ is true, Conjecture 2 is true for $\Delta_3(x)$ (see Remark 18 of \cite{CTZ2}). When $d=3,$ we can prove some partial results to Conjecture 2
(see \cite{CTZ2, CTZ4}).

{\bf Conjecture 3.} {\it Suppose $\ell\geq 2$ is a fixed integer. Suppose $\mathfrak{L}_f(s)\in \mathcal{S}$ is a function of degree $d\geq 2$, and $\mathfrak{L}(s)\in S$ is a function such that the asymptotic formula
\begin{equation}\label{wu.3}
\int_1^T|\mathfrak{L}(1/2+it)|^2dt=TP(\log T)+O(T^{1-\eta})
\end{equation}
for some $\eta>0,$ where $P(t)$ is a polynomial of degree $m.$

Then there exist constants
$\delta_1>0$ and $\delta_2>0$ such that the estimate
\begin{equation}\label{C-4}
\int_1^TE_f(x)|\mathfrak{L}(1/2+ix)|^2dx\ll  T^{1+\frac{ (d-1)}{2d}-\delta_1+\varepsilon}
\end{equation}
and the asymptotic formula
\begin{equation}\label{wu.5}
\int_1^TE_f^\ell(x)|\mathfrak{L}(1/2+ix)|^2dx=C_0C_{f,d}^{(\ell)}T^{1+\frac{\ell (d-1)}{2d}}Q(\log T)+O(T^{1+\frac{\ell (d-1)}{2d}-\delta_2+\varepsilon})
\end{equation}
hold, where $C_0$ is a constant independent of $f$, $C_{f,d}^{(\ell)}$ was defined in (\ref{C-2})     and  $Q(t)$ is a polynomial in $t$ of degree $m.$}

{\bf Remark 3.} From (\ref{er.16}) and  (\ref{er.18}) we see that if $\ell\leq 8$ then Conjecture 3 is true for $\mathfrak{L}_f(s)=\zeta^2(s), \  \mathfrak{L}(s)=\zeta(s)$. From the approach of Ivi\'c and Zhai \cite{IZ2} we see that Conjecture 3  is true for
any $\mathfrak{L}_f(s)\in \mathcal{S}$ of degree $2$ and  $\mathfrak{L}(s)=\zeta(s)$
when $\ell$ is small (say, $\ell\leq 8$). The conjecture (\ref{er.19}) is a special case of (\ref{wu.5}) in Conjecture 3.

The following question is an analogue of Problem 1 in Section 2.

{\bf Problem 2.} {\it Let   $\mathfrak{L}_{f}(s)\in \mathcal{S}$ be  a function of degree $d\geq 2.$  Do there exist intervals $I = [X, X + X^\beta]$ with $\beta>0$ such that the estimate
\begin{equation}\label{tsang-problem2}
 \int_{I}|E_{f} (x) | dx\ll
  T^{\beta+ \frac{ (d-1)}{2d}-\delta}
 \end{equation}holds
for some small positive $\delta$?}

\subsection{Conjectures on hybrid Moments of several error terms }\

{\bf Conjecture 4.} {\it Let $k\geq 2$   and $2\leq d_1<d_2<\cdots <d_k$ be fixed integers. For each $1\leq j\leq k,$ suppose $\mathfrak{L}_{f_j}(s)\in \mathcal{S}$ is  a function of degree $d_j\geq 2.$
Suppose $\ell_j\geq 1\ (j=1,2,\cdots,k)$  are  fixed integers and let
$\ell=\min(\ell_1,\ell_2,\cdots, \ell_k).$

If $\ell=1$, then there exists a  constant
$\delta_1>0$   such that the estimate
\begin{equation}\label{wu.7}
\int_1^TE_{f_1}^{\ell_1}(x)E_{f_2}^{\ell_2}(x) \cdots
E_{f_k}^{\ell_k}(x)dx\ll  T^{1+\sum_{j=1}^k\frac{\ell_j (d_j-1)}{2d_j}-\delta_1+\varepsilon}
\end{equation}
holds.

If $\ell>1$, then there exists a  constant   $\delta_2>0$ such that
 the asymptotic formula
\begin{equation}\label{wu.8}
\int_1^T E_{f_1}^{\ell_1}(x)E_{f_2}^{\ell_2}(x) \cdots
E_{f_k}^{\ell_k}(x)dx=C   T^{1+\sum_{j=1}^k\frac{\ell_j (d_j-1)}{2d_j}}
+O(T^{1+\sum_{j=1}^k\frac{\ell_j (d_j-1)}{2d_j}-\delta_2+\varepsilon})
\end{equation}
holds, where $$C=C_0\prod_{j=1}^kC_{f_j,d_j}^{(\ell_j)}\not= 0$$ and $C_0$ is a constant independent of $f_1, \cdots, f_k.$}

{\bf Remark 4.} We see that   results of Subsection 4.1 give some partial
answers to Conjecture 4. Note that all results in Subsection 4.1 remain valid if we replace $\zeta^2(s)$ by any function of degree two in $\mathcal{S}.$

{\bf Conjecture 5.} {\it Let $k\geq 2$   and $2\leq d_1<d_2<\cdots <d_k$ be fixed integers. For each $1\leq j\leq k,$ suppose $\mathfrak{L}_{f_j}(s)\in \mathcal{S}_{real}$ is  a function of degree $d_j\geq 2.$ Suppose $T\geq 10$ is a large real number.

(1) There exist  a positive constant $\delta_1<1$ and two points
 $x_\ell\in [T, T+  T^{\delta_1}](\ell=1,2)$ such that
 $$\prod_{j=1}^kE_{f_j}(x_1) >c T^{\sum_{j=1}^k\frac{d_j-1}{2d_j}}, \ \ \ \prod_{j=1}^kE_{f_j}(x_2)<-c T^{\sum_{j=1}^k\frac{d_j-1}{2d_j}} ,$$  where   $c $ is a   positive constant. We conjecture that any $\delta_1>1-1/d_k $ is OK.

 (2) There is a positive constant $0<\delta_2<1$ such that the following property holds:   in the interval $[T, 2T]$, there are at least $ T^{1-\delta_2}$ disjoint subintervals of length $\gg T^{\delta_2}$   such that
$$\pm\prod_{j=1}^kE_{f_j}(x) >c^{*} T^{\sum_{j=1}^k\frac{d_j-1}{2d_j}},$$  whenever $x$ lies in any of
these subintervals, where $c^{*}>0$ is a positive constant. We conjecture that any $\delta_2<1-1/d_1$ is OK. }

{\bf Remark 5.} From Theorem 1 and Theorem 2 we see that Conjecture 5 is true
for $k=2, \  \mathfrak{L}_{f_1}(s)=\zeta^2(s), \ \mathfrak{L}_{f_2}(s)=\zeta^3(s),\ \ \delta_1=139/144+\varepsilon,\ \delta_2=13/48-\varepsilon$(or even $\delta_2=1/2-\varepsilon$ under some strong conditions).
Actually we can prove that if $\mathfrak{L}_{f_1}(s)\in \mathcal{S}_{real}$ is any function of degree $2$, then
Conjecture 5 is true for $k=2, \  \mathfrak{L}_{f_1}(s),\  \mathfrak{L}_{f_2}(s)=\zeta^3(s)$      with suitable $\delta_\ell>0(\ell=1,2)$.

{\bf Conjecture 6.} {\it Let $k\geq 2$   and $2\leq d_1<d_2<\cdots <d_k$ be fixed integers. For each $1\leq j\leq k,$ suppose $\mathfrak{L}_{f_j}(s)\in \mathcal{S}$ is  a function of degree $d_j\geq 2.$ Suppose $\ell_j\geq 1\ (j=1,2,\cdots,k)$  are  fixed integers and $\ell=\min(\ell_1,\ell_2,\cdots, \ell_k).$ Suppose $\mathfrak{L}(s)\in S$ is a function such that the asymptotic formula (\ref{wu.3}) holds.

If $\ell=1$, then there exists a  constant
$\delta_1>0$   such that the estimate
\begin{equation}
\int_1^TE_{f_1}^{\ell_1}(x)E_{f_2}^{\ell_2}(x) \cdots
E_{f_k}^{\ell_k}(x)|\mathfrak{L}(1/2+ix)|^2dx\ll  T^{1+\sum_{j=1}^k\frac{\ell_j (d_j-1)}{2d_j}-\delta_1+\varepsilon}
\end{equation}
holds.

If $\ell>1$, then there exists a  constant   $\delta_2>0$ such that
 the asymptotic formula
\begin{eqnarray}
&&\ \ \ \ \ \int_1^T E_{f_1}^{\ell_1}(x)E_{f_2}^{\ell_2}(x) \cdots
E_{f_k}^{\ell_k}(x)|\mathfrak{L}(1/2+ix)|^2dx\\
&&=C_*   T^{1+\sum_{j=1}^k\frac{\ell_j (d_j-1)}{2d_j}}\mathcal{Q}(\log T)
+O(T^{1+\sum_{j=1}^k\frac{\ell_j (d_j-1)}{2d_j}-\delta_2+\varepsilon})\nonumber
\end{eqnarray}
holds, where $$C_*=C_1\prod_{j=1}^kC_{f_j,d_j}^{(\ell_j)}\not= 0$$ with $C_1$   a constant independent of $f_1, \cdots, f_k,$ and $\mathcal{Q}(t)$ is a polynomial in $t$ of degree $m.$}

The following Problem 3 is an analogue and generalization of Problem 2.

{\bf Problem 3.} {\it Let $k\geq 2$ be a fixed integer and $2\leq d_1<d_2<\cdots <d_k$ be fixed   integers. For each $1\leq j\leq k,$ suppose $\mathfrak{L}_{f_j}(s)\in \mathcal{S}$ is  a function of degree $d_j\geq 2.$  Do there exist intervals $I = [X, X + X^\beta]$ with $\beta>0$ such that the estimate
\begin{equation}\label{tsang-problem3}
 \int_{I}|E_{f_1} (x)E_{f_2} (x) \cdots E_{f_k} (x)| dx\ll
  T^{\beta+\sum_{j=1}^k\frac{ (d_j-1)}{2d_j}-\delta}
 \end{equation}holds
for some small positive $\delta$?}

{\bf Acknowledgement} The author deeply thanks the anonymous referee for valuable suggestions.

{\bf Conflict of Interest} The author declares that there is no conflict of interests regarding the publication of this paper.

\end{document}